\documentclass[titlepage,12pt]{article}

\newcommand{\calR}{{\cal{R}}}

\newcommand{\HH}{{\Bbb{H}}}

\renewcommand{\H}{\HH}
\newcommand{\calP}{{\cal{P}}}
\newcommand{\calL}{{\cal{L}}}
\newcommand{\syst}{{\rm{syst}}}
\newcommand{\U}{{\cal{U}}}
\newcommand{\V}{{\cal{V}}}

\newtheorem{Def}{Definition}[section]

\newtheorem{Th}{Theorem}[section]
\newenvironment{Thnmb}{\vskip 6pt \noindent
{\bf{Theorem  \nmb}} \begin{it}}{\end{it} \vskip 6pt}

\newenvironment{Ack}{\vskip 6pt \noindent{\bf
Acknowledgements:} \quad}{\vskip 6pt}

\newtheorem{Cor}{Corollary}[section]   
\newtheorem{Lm}{Lemma}[section] 
\newcommand{\Pf}{\noindent
{\bf{Proof}}:\quad}

\newcommand{\vol}{{\hbox{ \rm vol}}}

\newcommand{\Z}{{\Bbb{Z}}}

\newcommand{\calF}{{\cal{F}}}
\newcommand{\calO}{{\cal{O}}}

\newcommand{\const}{{\hbox{\rm const}}}
\newcommand{\tr}{{\hbox{tr}}}

\newcommand{\area}{{\hbox{area}}}

\newcommand{\eps}{\varepsilon}

\newcommand{\diam}{{\hbox{\rm diam}}}
\newcommand{\length}{{\hbox{length}}}

\newcommand{\C}{{\Bbb{C}}}

\newcommand{\btbt}{\left( \begin{array}{cc}}
\newcommand{\etbt}{ \end{array}\right)}
\newcommand{\bthbth}{\left( \begin{array}{ccc}}
\newcommand{\ethbth}{ \end{array}\right)}

\newcommand{\bcol}{\left(\begin{array}{c}}
\newcommand{\ecol}{\end{array}\right)}

\newcommand{\bthree}{\left( \begin{array}{ccc}}

\newcommand{\genus}{{\hbox{\rm genus}}}

\usepackage{amsfonts}
\usepackage{latexsym}
\usepackage[dvips]{graphicx}
\usepackage{epsf}
\begin{document}

\title{Random Construction of Riemann Surfaces}

\author{Robert Brooks\thanks{Partially supported by the Israel Science
Foundation, the Fund for the Support of Research and the Fund for the Support
of Sponsored Research at the Technion, and the M. and M. L. Bank
Research Fund}\\ Department of Mathematics\\
Technion -- Israel Institute of Technology\\ Haifa, Israel\\
rbrooks@tx.technion.ac.il  \and Eran
Makover\thanks{Partially supported by the NSF Grant DMS-0072534} \\ Department of Mathematics\\
Dartmouth College\\ Hanover, New Hampshire 03755\\ Eran.Makover@dartmouth.edu}

\date{August, 2000}

\maketitle

\newcommand{\calFns}{{\cal{F}}_n^*}
\newcommand{\Prob}{{\hbox{Prob}}}
\renewcommand{\Bbb}{\mathbb}
In this paper, we address the following question: What does a typical
compact Riemann surface of large genus look like geometrically? 

By a Riemann surface, we mean an oriented surface with a complete,
finite-area  metric of
constant curvature {-1}. 

In the standard geometric picture of Riemann surfaces via
Fenchel-Nielsen coordinates, it is difficult to keep track of global
geometric quantities such as the first eigenvalue of the Laplacian,
the injectivity radius, and the diameter.  Here, we present a model
for looking at Riemann surfaces based on $3$-regular graphs, with
which it is easier to control the global geometry.

The idea
of using $3$-regular graphs to study the first eigenvalue of Riemann
surfaces originated in Buser (\cite{Bu1}, \cite{Bu2}), who associated
cubic graphs to Riemann surfaces as a tool for comparing  the
spectral geometry of surfaces with the spectral geometry of graphs.  We
introduce a somewhat different method, which associates to each
$3$-regular  graph with
an orientation a finite area Riemann surface.

To be more specific, if $\Gamma$ is a finite $3$-regular graph, an
orientation $\calO$ on $\Gamma$ is a function which assigns to
each vertex $v$ of $\Gamma$ a cyclic ordering of the edges emanating
from $v$. In \S 2 below, we will show how, given a pair 
$(\Gamma, \calO)$, we may associate to $(\Gamma, \calO)$ two Riemann
surfaces $S^O(\Gamma, \calO)$ and $S^C(\Gamma, \calO)$. $S^O(\Gamma,
\calO)$ is constructed by associating an ideal hyperbolic triangle to each
vertex of $\Gamma$, and gluing sides together according to the edges
of the graph $\Gamma$ and the orientation $\calO$. It is a finite-area
Riemann surface with cusps. 

The surface $S^C(\Gamma, \calO)$ is then the conformal
compactification of $S^O(\Gamma, \calO)$.

It follows from a theorem of Belyi \cite{Bel} (see \cite{JS} for a
 discussion of Belyi's Theorem),
 that the surfaces
$S^C(\Gamma, \calO)$ are dense in the space of all Riemann 
surfaces. Thus, the
process of randomly selecting a Riemann surface can be modeled on the
process of  picking a
finite $3$-regular graph with orientation at random.

Since the pair $(\Gamma, \calO)$ gives a description of
 $S^O(\Gamma, \calO)$ as an
orbifold covering of $\H^2/PSL(2, \Z)$, one can give a qualitative
description of the global geometry of $S^O(\Gamma, \Z)$ by a
corresponding description of the pair $(\Gamma, \calO)$. Thus, for
example, the first eigenvalue of $S^O(\Gamma, \calO)$ will be large if
and only if the first eigenvalue of $\Gamma$ is large (\cite{SGTC} and
 \cite{VD}). 

It is our observation here, building on the work of \cite{PS}, that
the same will be true of $S^C(\Gamma, \calO)$, provided that
$S^O(\Gamma, \calO)$ satisfies a ``large cusps'' condition, to be
described in \S 1 below. This has a purely combinatorial interpretation in
terms of the pair $(\Gamma, \calO)$, and so can be analyzed with
relative ease.

For $n$  a positive integer, let $\calFns$ denote the finite set of
pairs $(\Gamma, \calO)$, where $\Gamma$ is a $3$-regular graph on $2n$
vertices. We will endow $\calFns$ with a probability measure
introduced and studied by Bollob\'as \cite{Bol1}, \cite{Bol2}, which we
review in \S 3 below. 

If $Q$ is
a property of $3$-regular graphs with orientation, denote by
$\Prob_n[Q]$ the probability that a pair $(\Gamma, \calO)$ picked from
$\calFns$ has property $Q$.

Our main technical result, shown in \S 4 below, is:

\begin{Th} \label{cusps}

As $n \to \infty$, 
$$\Prob_n [ S^O(\Gamma, \calO)\ {\hbox{satisfies the large cusps
condition}}\ ] \to 1.$$

\end{Th}

We will use Theorem \ref{cusps} in order to study geometric properties
of the surfaces $S^C(\Gamma, \calO).$ To that end, we define the
Cheeger constant $h(S)$ of a 
Riemann surface $S$ by the formula

$$h(S) = \inf_C {{\length(C)}\over{\min[\area(A), \area(B)]}},$$
where $C$ runs over (possibly disconnected) closed curves on $S$ which
divide $S$ into two parts $A$ and $B$.
 
It will then follow from Theorem \ref{cusps} that:
 
\begin{Th}\label{typical} There exist constants $C_1$, $C_2$, $C_3$,
and $C_4$ such that, as $n \to \infty$:

\begin{description}
\item{(a)} The first eigenvalue $\lambda_1(S^C(\Gamma, \calO))$
satisfies
 $$\Prob_n[ \lambda_1(S^C(\Gamma, \calO)) \ge C_1] \to 1.$$
 
\item{(b)}The Cheeger constant $h(S^C(\Gamma, \calO))$ satisfies 
 $$\Prob_n [h(S^C(\Gamma,\calO)) \ge C_2] \to 1.$$ 

\item{(c)} The shortest geodesic $\syst(S^C(\Gamma, \calO))$ satisfies
$$Prob_n[\syst(S^C(\Gamma, \calO)) \ge C_3] \to 1.$$

\item{(d)} The diameter $\diam(S^C(\Gamma, \calO))$ satisfies
$$\Prob_n[\diam(S^C(\Gamma, \calO)) \le C_4
\log(\genus(S^C(\Gamma, \calO)))] \to 1.$$
\end{description}
\end{Th}

Of these properties, (a) follows from (b) by Cheeger's inequality
\cite{Ch}, while (d) also follows from (b) and (c) and the following
well-known argument: if $M$ is a manifold and $B(r_0)$ is the infimum
over all points of $M$ of the volume of a ball of radius $r_0$, then
$$\diam(M) \le 2 [ r_0 + {{1}\over{h(M)}} \log( {{\vol(M)}\over{2
B(r_0)}})].$$ 

 Property
(b) for the surfaces $S^O(\Gamma, \calO)$ will  follow from the
corresponding result on graphs \cite{Bol2} together with \cite{SGTC}
and \cite{VD}, while the passage from the
 surfaces $S^O(\Gamma, \calO)$ to the surfaces $S^C(\Gamma, \calO)$
will follow from Theorem \ref{cusps}. For properties (c) and (d), the
translation from graphs to surfaces 
is not as simple, but the idea is similar.

The forms of (a), (b), and (d) are sharp, up to constants.  Regarding
(a), it follows from Cheng's Theorem \cite{Chng} that a Riemann
surface $S$ must have $\lambda_1(S) \le 1/4 + \eps$ for some $\eps \to
0$ as ${\hbox{genus}}(S) \to \infty$. The upper bound $h(S) \le 1 +
\eps$, is well-known, and follows from a similar argument.   The
estimate  $\diam(S) \ge (\const)
\log({\hbox{genus}}(S))$ follows from area considerations and
Gauss-Bonnet.  

The estimate in (c) is certainly not
optimal, as there are Riemann surfaces whose injectivity radius grows
like $(\const)[\log({\hbox{genus(S))}}]$. Indeed, this occurs for the
Platonic surfaces \cite{PS}, and also for congruence
coverings of compact arithmetically defined surfaces.  It follows from our
analysis that, for a given constant $C_5$, there is a positive constant $C_6$
such that
$$Prob_n[ \syst(S^C(\Gamma, \calO)) \ge C_5] \to C_6.$$
Thus, probability of selecting a surface having injectivity radius at
least a given large number is asymptotically positive, but certainly
not asymptotically 1.

In the language of \cite{SGNT}, Theorem \ref{typical} shows that, with
 probability $\to 1$, a typical Riemann surface is short, with a large
first eigenvalue, but not necessarily fat. 

The results that we use from \cite{PS} are qualitative rather than
quantitative. However, they have been put in quantitative form in
\cite{M}. In particular, it follows from \cite{M} that whenever in the
following the
condition of ``cusps of length $\ge L$'' is used, we may take $L=7$.

The results of Theorem \ref{typical} were announced in \cite{BM},
under the weaker conclusion that properties (a)-(d) occur with
positive  probabilty, rather than probability $\to 1$ as $n \to
\infty$.

\begin{Ack} The first author acknowledges with gratitude the
hospitality of the \'Ecole Normale Superiere of Lyon, where much of
the writing was carried out, and also the hospitality of Dartmouth
College. 

\end{Ack}

\section{Compactification of Riemann Surfaces}

In this section, we  review the connection between a finite-area
Riemann surface and its conformal compactification.

Let $S^0$ be a
Riemann surface with a complete finite area metric of curvature $-1$. 
Then $S^0$ has finitely many cusps neighborhoods $C^1, \ldots, C^k$,
such that, for each $C^i$ there is an isometry $$f_i: C^i \to
C^{y_i}$$ for some $y_i$, where $C^{y_i}$ is the space $$C^{y_i} = \{
z \in \C : \Im(z) \ge {{1}\over{y_i}}\}/ (z \sim z
+1),$$
endowed with the hyperbolic metric 
$$ds^2 = {{1}\over{y^2}}\left[ dx^2 + dy^2 \right].$$

The curve $h_i = f_i^{-1}(z: \Im(z) = {{1}\over{y_i}})$ on $S^O$ is a
closed horocycle on $C^i$ whose length is $y_0$. The length of the
largest closed horocycle on the cusp $C^i$ is a measure of how large
the cusp is.

\begin{figure}[hbt] 
\leavevmode
\vskip 12pt
\centerline{\epsfysize 2.5in \epsfbox{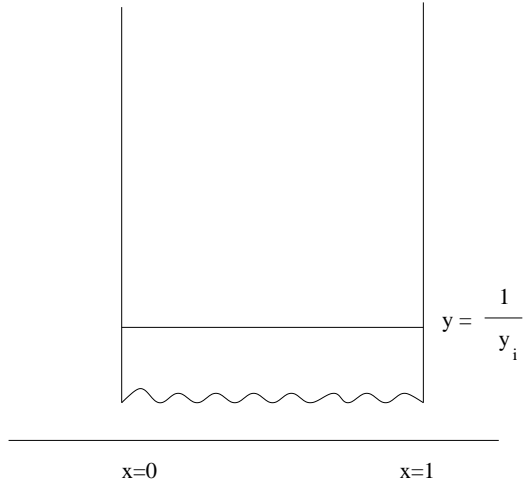}}
\caption{The length of a horocycle}\label{horoc}
\end{figure}

\begin{Def} $S^O$ has cusps of length $\ge L$ if we may choose all the
$C^i$'s disjoint, with $y_i \ge L$ for all $i= 1, \ldots, k$.

\end{Def}

Given a finite-area Riemann surface $S^O$, there is a unique compact
Riemann surface 
$S^C$ and finitely many points $p_1, \ldots, p_k$ of $S^C$ such that $S^O$ is
conformally equivalent to $S^C - \{p_1, \ldots, p_k\}$. $S^C$ may be
constructed from $S^O$ by observing that each cusp neighborhood $C^i$
is conformally equivalent to a punctured disk. One may fill in this
puncture conformally and then reglue the disk to obtain the conformal 
structure on the closed surface $S^C$. By the Uniformization Theorem,
there is a unique constant 
curvature metric which agrees with this conformal structure.

It is natural to raise the question of the relationship between the
constant curvature metric on $S^O$ and the constant curvature metric
on $S^C$.
In general, the relationship need not be close. For instance, the
surface $S^C$ need not carry a hyperbolic metric even when $S^O$
carries one. However,
it is shown in Theorem 2.1 of \cite{PS} that, if $S^O$ has cusps of
length $ \ge L$ for a suitably large $L$, then $S^C$ will carry a
hyperbolic metric, and indeed  this metric will be very closely
related to the hyperbolic metric on $S^O$. 

More precisely, we have:

\begin{Th}[{\cite{PS}}] For every $\eps$, there exists numbers $L, r,$ 
and $y$ such
that, if the cusps of $S^O$ have length $\ge L$, then, outside the
union of cusp neighborhoods $\U= \cup_{i=1}^k f_i^{-1}(C^{y})\subset 
S^{O}$ of 
the cusps $C^{i}$, and $\V= \cup_{i=1}^k B_{r}(p_i)\subset S^{C}$, the metrics
$ds_C^2$ and $ds_O^2$ satisfy
$${{1}\over{(1+ \eps)}} ds_O^2 \le ds_C^2 \le (1 + \eps)ds^2_O.$$
\end{Th}

The proof of this theorem is based on the Ahlfors Schwarz Lemma \cite{A}. 
The idea of the proof
 is to build on the compact surface an intermediate metric
$ds^2_{{\hbox{int}}}$ 
 with  curvature close to the curvature of the 
metric on the open surface, and  to 
use the Ahlfors-Schwarz Lemma  to compare this metric to the constant
curvature metric. The large cusps condition enters precisely here,
by giving the metric $ds^2_{{\hbox{int}}}$ sufficient time to evolve
from the hyperbolic metric on the ball to the hyperbolic metric on the
cusp, while keeping curvature close to constant.

It was shown in \cite{PS} that this result may
be employed to show that, under the assumption of large cusps, the
surfaces $S^O$ and $S^C$ share a number of global geometric
properties.

\begin{Th}[{\cite{PS}}]\label{ps} For every $\eps$, there exists an $L$ such
that, if $S^O$ has cusps of length $\ge L$, then

\begin{description} 

\item{(a)} the Cheeger constants $h(S^O)$ and $h(S^C)$ satisfy
$${{1}\over{(1 + \eps)}} h(S^O) \le h(S^C) \le (1 + \eps) h(S^O).$$

\item{(b)} the shortest closed geodesics $\syst(S^O)$ and $\syst(S^C)$
satisfy
$${{1}\over{(1 + \eps)}} \syst(S^O) \le \syst(S^C).$$

\end{description}
\end{Th}
We do not obtain an inequality of the form
$$\syst(S^C) \le (\const) \syst(S^O),$$
in (b), because it may happen that the shortest closed geodesic on $S^O$
becomes homotopically trivial on $S^C$.

\section{Surfaces and $3$-Regular Graphs}

Let $\Gamma$ be a finite $3$-regular graph. An orientation $\calO$ on
the graph is the assignment, for each vertex $v$ of $\Gamma$, of a
cyclic ordering of the three edges emanating from $v$. If $\Gamma$ has
$2n$ vertices, then clearly there are $2^{2n}$ orientations on $\Gamma$.

We may think of an orientation on a $3$-regular graph in the following
way: suppose one were to walk along the graph. Then, when one
approaches a vertex, the orientation allows one to distinguish between
a left-hand turn and a right-hand turn at the vertex. Thus, any path
on $\Gamma$ beginning at a vertex $v_0$ may be described by picking an
initial direction and a series of $L$'s (signalling a left-hand turn)
and $R$'s (signalling a right-hand turn).

To a pair $(\Gamma, \calO)$ we will associate two Riemann surfaces
$S^O(\Gamma, \calO)$ and $S^C(\Gamma, \calO)$ as follows: We begin by
considering the ideal hyperbolic triangle $T$ with vertices $0, 1$,
and $\infty$ shown in Figure \ref{triangle} below. The solid lines in
Figure \ref{triangle} are geodesics joining the points $i, i+1$, and
${{i+1}\over{2}}$ with the point ${{1 + i \sqrt{3}}\over{2}}$, while
the dotted lines are  horocycles of length 1 joining pairs of points
from the set $\{i, i+1, {{i+1}\over{2}}\}$. We may think of these
points as ``midpoints'' of the corresponding sides, even though they
are of infinite length. We may also think of the three solid lines as
segments of a graph surrounding a vertex. We  then give them the
cyclic ordering $(i, i+1, {{i+1}\over{2}})$.

\begin{figure}[!h]
\begin{center}

\includegraphics[angle=0,scale=0.4]{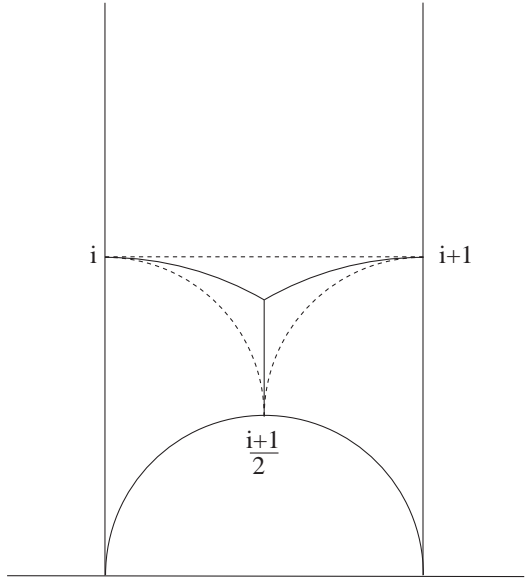}

\caption{The marked ideal triangle $T$}\label{triangle}
\end{center}
\end{figure}

We may now construct the surface $S^O(\Gamma, \calO)$ from $(\Gamma,
\calO)$ by placing on each vertex $v$ of $\Gamma$ a copy of $T$, so
that the cyclic ordering of the segments in $T$ agrees with the
orientation at the vertex $v$ in $\Gamma$. If two vertices of $\Gamma$
are joined by an edge, we glue the two copies of $T$ along the
corresponding sides subject to the following two conditions:

\begin{description}

\item{(a)} the midpoints of the two sides are glued together,

\item{\ { }} and

\item{(b)} the gluing preserves the orientation of the two copies of
$T$.

\end{description}

The conditions (a) and (b) determine the gluing uniquely. It is
easily seen that the surface $S^O(\Gamma, \calO)$ is a complete
Riemann surface with finite area equal to $2\pi n$, where $2n$ is the
number of vertices of $\Gamma$.

The surface $S^C(\Gamma, \calO)$ is then the conformal
compactification of $S^O(\Gamma, \calO)$.

In the remainder of this section, we will describe how to read off
many geometric properties of the surfaces $S^O(\Gamma, \calO)$ and
$S^C(\Gamma, \calO)$ from the combinatorics of the pair $(\Gamma,
\calO)$.

We begin with the observation that the topology of $S^O$ is easy to
reconstruct from $(\Gamma, \calO)$.  We will need the following

\begin{Def} A left-hand-turn path on $(\Gamma, \calO)$
 is a closed path on $\Gamma$
such that, at each vertex, the path turns left in the orientation
$\calO$.

\end{Def}

Traveling on a path on $\Gamma$ which always turns left describes a
path on $S^O(\Gamma, \calO)$ which travels around a cusp. Indeed, if
we set $l= l(\Gamma, \calO)$ to be the number of disjoint
left-hand-turn paths, then the topology of $S^O(\Gamma, \calO)$ is
easily describable in terms of $l$ and $n$. Indeed, the graph $\Gamma$
divides $S^O(\Gamma, \calO)$ into $l$ regions, each bordered by a
left-hand-turn path and containing one cusp in its interior.
From this, we can immediately read off the signature of $S^O(\Gamma,
\calO)$ 
by the Euler
characteristic. The genus of $S^O(\Gamma, \calO)$ is given by
$$\genus = 1 + {{n-l}\over{2}},$$
and the number of cusps is $l$.

Note that the usual orientation on the $3$-regular graph which is the
1-skeleton of the cube contains six left-hand-turn paths, giving that
the associated surface is a sphere with six punctures, while a choice
of a different orientation on this graph can have either two, four, or
six left-hand-turn paths, so that the associated surface can have
genus 0, 1, or 2. Thus, the topology of $S^O(\Gamma, \calO)$ is
heavily dependent on the choice of $\calO$.

The geometry of the cusps can also be read off from $(\Gamma,
\calO)$. To that end, we observe that the horocycles on the various
copies of $T$ fit together to form a system of disjoint closed
horocycles about the cusps of $S^O(\Gamma, \calO)$. We call this
system of horocycles the {\em canonical horocycles} of $S^O(\Gamma,
\calO)$. The length of each closed horocycle in this set is precisely
the length of the corresponding left-hand-turn path, since the length
of the horocycle joining $i$ to $i+1$ has length 1. Thus, $S^O(\Gamma,
\calO)$ has cusps of length $\ge L$ if the length of any left-hand
turn path on $(\Gamma, \calO)$ is at least $L$.

The converse to this is not true, as it is possible to choose a system
of horocycles other than the canonical horocycles such that the length
of the shortest horocycle for the new system is larger than the length
of the shortest canonical horocycle. We will return to this idea in \S
4 below.

The Cheeger constant $h(S^O(\Gamma, \calO))$ can be estimated in terms
of the graph $(\Gamma, \calO)$ as well.
Recall that, by analogy with the Cheeger constant of a manifold,
 the Cheeger constant $h(\Gamma)$ of a graph $\Gamma$ is
given by
$$h(\Gamma) = \inf_E {{\#(E)}\over{\min(\# (A), \#(B))}},$$
where $E$ is a collection of edges such that $\Gamma - E$ disconnects
into two components $A$ and $B$, and $\#(A)$ (resp. $\#(B)$) is the
number of vertices in $A$ (resp. $B$).

Then we  have

\begin{Th}[\cite{SGTC}, {\cite{VD}}]\label{cheeger} 
There are positive constants
$C_1$ and $C_2$ such that
$$C_1 h(\Gamma) \le h(S^O(\Gamma, \calO)) \le C_2 h(\Gamma)$$
for all finite $3$-regular graphs $\Gamma$.
\end{Th}

In effect, the pair $(\Gamma, \calO)$ describes $S^O(\Gamma, \calO)$
as an orbifold covering space of the orbifold $\H^2/PSL(2,  \Z)$. The
behavior of the Cheeger constant of a finite covering of a compact
manifold in terms of the graph of a covering is described in
\cite{SGTC}. In the present case, the base manifold is not compact,
but rather a finite-area Riemann surface (with singularities). The
additional complication which this difficulty presents is solved in
\cite{VD}. 

This gives Theorem \ref{cheeger}.

Note in particular that the quantity
$h(\Gamma)$  of
Theorem \ref{cheeger} depends only on $\Gamma$ and not on $\calO$.

The geodesics of $S^O(\Gamma, \calO)$ are also describable in terms of
$(\Gamma, \calO)$. To explain this, let $\calL$ and $\calR$ denote the
matrices 
$$\calL = \btbt 1&1\\0&1 \etbt \quad \calR = \btbt 1 & 0 \\ 1 & 1
\etbt.$$
\renewcommand{\l}{\length}

A closed path $\calP$ of length $k$ on the graph may be described by
starting at a midpoint of an edge, and then giving a sequence $(w_1,
\ldots, w_k),$ where each $w_i$ is either $l$ or $r$, signifying a left
or right turn at the upcoming vertex. We then consider the matrix
$$M_{\calP} = W_1 \ldots W_k,$$
where $W_j = \calL$ if $w_j = l$ and $W_j = \calR$ if $w_j=r$. 

The closed path $\calP$ on $\Gamma$ is then homotopic to a closed
geodesic $\gamma(\calP)$ on $S^O(\Gamma, \calO)$ whose length
$\l(\gamma(\calP))$ is given by
$$2\cosh({{\l(\gamma(\calP))}\over{2}})= \tr(M_{\calP}).$$

Note that $\l(\gamma(\calP))$ depends very strongly on $\calO$. Indeed,
if the path $\calP$ contains only left-hand turns, then
$\l(\gamma(\calP))=0$, and if $\gamma(\calP)$ is a path of length $r$
containing  precisely
one right-hand turn, then 
$$\l(\gamma(\calP))= 2  \log\left( {{(1+r) + \sqrt{(1+r)^2 -
4}}\over{2}}\right) \sim 2\log(1 +r),$$
and hence grows linearly in $\log(r)$.
On the other hand, if the path $\calP$ of length $r$ 
consists of alternating left-
and right-hand turns, then
$$\l(\gamma(\calP)) = r \log( {{3 + \sqrt{5}}\over{2}}),$$
which is linear in $r$.

We now consider the description of $S^C(\Gamma,\calO)$ in terms of
$(\Gamma,\calO)$. We will carry out this description under the
assumption that the cusps of $S^O(\Gamma, \calO)$ are large. 

\begin{Th}\label{compact} Assume that the cusps of $S^O(\Gamma,
\calO)$ have length $\ge L= L(\eps)$. Then there exist constants $C_1,
C_2, C_3, C_4$, and $C_5$  depending only on $L$ such that: 

\begin{description}

\item{(a)} The Cheeger constant $h(S^C(\Gamma, \calO))$ satisfies
$$C_1 h(\Gamma) \le h(S^C(\Gamma, \calO)) \le C_2 h(\Gamma).$$

\item{(b)} The shortest closed geodesic $\syst(S^C(\Gamma, \calO))$ satisfies
$$ \syst(S^C(\Gamma, \calO)) \ge C_3 \log( 1 + \syst(\Gamma)) \ge C_4,$$
where $\syst(\Gamma)$ is the girth of the graph $\Gamma$.

\item{(c)} The genus of $S^C(\Gamma, \calO)$ satisfies
 $$\genus(S^C(\Gamma, \calO)) \ge C_5 \#(\Gamma).$$
\end{description}

\end{Th}

\Pf (a) follows from Theorem 1.2 (a) and Theorem 2.1.

(b) follows from Theorem 1.2 (b) and the calculation of lengths of
geodesics on $S^O(\Gamma, \calO)$. 

(c) follows from the formula for the genus of $S^O(\Gamma, \calO)$,
which is also the genus of $S^C(\Gamma, \calO)$, together with the
simple observation that if each cusp in $S^O(\Gamma, \calO)$ is
bounded by a horocycle of length at least $L$, then the number of such
cusps is bounded by ${{1}\over{L}}[\area(S^O(\Gamma, \calO))]$, since
$L$  is the area inside a horocycle of length $L$.

%\section{A Model for Selecting Random $3$-Regular Graphs}
\section{The Bollob\'as Model}

In this section we discuss a model, due to Bollob\'as, for studying the
process of randomly selecting a $3$-regular graph.

The problem of putting a probability measure on the set of $3$-regular
graphs on $2n$ vertices would not appear at first sight to be
difficult, since this is a finite set. It has, however, proven
problematic to find as model which is amenable to meaningful
calculation, and a number of different models have been proposed and
studied, each with its own benefis and drawbacks. See Janson
\cite{Jan} for a  discussion and comparison of the different models.

 We will use a model  introduced by Bollob\'as
(\cite{Bol1}, \cite{Bol2}).
Bollob\'as considered the problem for $k$-regular graphs, $k$ arbitrary,
 but we will need only
the case $k=3$. This model has the advantage that calculations of an
asymptotic character (as $n \to \infty$) can be carried out with
relative ease.

For each $n$, let $\calF_n$ denote the finite set of $3$-regular graphs on
$2n$ vertices, and let $\calF_n^*$ denote the set of $3$-regular graphs
with orientation. We put a probability measure on $\calF_n$ and
$\calF_n^*$ in the following way: we consider a hat with $6n$ balls, each
ball labeled with an integer between $1$ and $2n$, each integer occurring
three times. We then build a graph at random by selecting pairs of balls
from the hat, without replacement. If at step $i$ the numbers $b_i$ and
$c_i$ are selected, we add to the graph an edge joining $b_i$ and $c_i$.
An orientation at vertex $b$ is determined by the cyclic order in which
edges are added to the vertex $b$.

We will need two results of Bollob\'as concerning this model.

The first result concerns the Cheeger constant of a graph:

\begin{Th}[{\cite{Bol2}}] There is a constant $C>0$ such that the
probability of a graph $\Gamma$ chosen randomly from $\calF_n$ having
Cheeger constant $h(\Gamma)$ greater than $C$ satisfies
$$\Prob_n [h(\Gamma) >C] \to 1\ {\hbox{as}}\ n \to \infty.$$
\end{Th}

Bollob\'as gives numerical estimates showing that $C> 2/11$.

To state the second result, we recall the notion of an asymptotic
Poisson distribution:

\begin{Def}
\begin{description}
\item{(a)} A random variable $X$ which takes values in the natural
numbers $\Z^+$  is a Poisson distribution with mean $\mu$ if
$$\Prob (X=k) = e^{-\mu} {{\mu^k}\over{k!}}.$$
The mean $\mu$ is the expected value of $X$.

\item{(b)} Let $\{X^n\}$ be a family of random variables on the
probability spaces $\{P_n\}$. 

The $\{X^n\}$ are asymptotic Poisson distributions as $n \to
\infty$ if there exists $\mu$ such that
$$\lim_{n \to \infty} \Prob( X^n = k) = e^{-\mu} {{\mu^k}\over{k!}}$$  
for all $k$.

\item{(c)} The families $\{X_i^n\}$ are asymptotically independent
Poisson distributions if, for each $i$, the random variables $X_i^n$
tend to a Poisson distribution $X_i$ as $n\to \infty$, and if the
variables $X_i$ are independent.
\end{description}
\end{Def}

A well-kmown example of an asymptotic Poisson distribution is
given by the hatcheck lady who returns hats in a random fashion to the
$n$ guests at a party. The random variable $X^n$ which is the number
of guests who receive the correct hat is asymptotically Poisson with
mean $1$ as $n \to \infty$.

Bollob\'as proves:

\begin{Th}[{\cite{Bol1}}] \label{poisson} Let $X_i$ denote the number of
closed paths in $\Gamma$ of length $i$. Then the random variables $X_i$ on
$\calF_n$ are asymptotically independent Poisson distributions with means
$$\lambda_i = {{2^i}\over{2i}}.$$

\end{Th}

In our case we have an additional structure on the graph- ``the
orientation''.  We are distinguishing between short paths that agree
with the orientation and those that do not.  To do so we look at all
the possible orientations on a given closed path of length $i$.  There
are $2^{i}$ possible orientations, but only two yield a left hand path.
Therefore we get:

\begin{Cor}\label{lcusps} Let $Y_i$ be the random variable on $\calF_n^*$
which associates to $(\Gamma, \calO)$ the number of left-hand turn paths
of length $i$. Then the $Y_i$ are asymptotically independent Poisson
distributions with means
$$\mu_i = {{1}\over{i}}.$$
\end{Cor}

Theorem \ref{poisson} and Corollary \ref{lcusps} imply that short
geodesics and small cusps will occur with positive probablity in the
surfaces $S^O(\Gamma, \calO)$, asymptotically as $n \to \infty$.
 One would expect on the grounds of
asymptotic independence that as $n \to \infty$, these phenomena appear
far apart. The following elementary lemma makes this expectation
precise (compare \cite{Bol1}, Theorem 32):

\begin{Lm} \label{dist} For fixed numbers $l_1, l_2,$ and $d$, let
$Q_n(l_1, l_2,d)$ denote the probability that a graph picked from
$\calF_n$ (resp.
$\calF_n^*$) has closed paths $\gamma_1$ and $\gamma_2$ of length $l_1$
and $l_2$ respectively, which are distance $d$ apart.

Then
$$Q_n(l_1, l_2, d) \to 0\ {\hbox{as}}\ n \to \infty.$$
\end{Lm}

\Pf We first observe that, since the statement is independent of the
orientation $\calO$, we may restrict our attention to picking from
$\calF_n$. We will show that, for every
$\eps$, for
$n$ sufficiently large, we have that
$$Q_n(l_1, l_2, d) < \eps.$$

Since the number of closed paths of length $l_1$ are asymptotically
Poisson distributed by Theorem \ref{poisson}, given $\eps_1$, we may find
$N(\eps_1)$ such that with probability $>1- \eps_1$, the number of closed
paths of length $l_1$ is less than $N(\eps_1).$

Now let $\gamma$ be a closed path in $\Gamma$ of length $l_1$. We
consider the $l_1\cdot 2^{d + [{{l_2}\over{2}}] -1}$ vertices which are
at distance at most $d+ [{{l_2}\over{2}}]$ from $\gamma$. When $n$ is
large compared to $l_1 \cdot 2^{d + [{{l_2}\over{2}}]}$, with probability
$\to 1$ as $n \to \infty$, no vertex in this set will have been selected
twice. This implies that there will be no closed path of length $l_2$ at
distance $d$ from $\gamma$. 

Applying this estimate to each of the $<
N(\eps_1)$ closed paths of length $l_1$ then gives the lemma.

\section{Large Cusps}

In this section, we complete the proof of Theorem \ref{cusps}. We will
reformulate it in the following way:

\def\nmb{\ref{cusps}A}
\begin{Thnmb}
Given $L$, as $n \to \infty$, we have
$$\Prob_n[ S^O(\Gamma, \calO)\ {\hbox{has cusps of length}}\ \ge L]
\to 1.$$
\end{Thnmb}

\Pf We begin the proof by calculating the probability that the
canonical horocycles of $S^O(\Gamma, \calO)$ all have length $\ge
L$. This is precisely the probability that all the random variables
$Y_i$ of Corollary \ref{lcusps} have the value $0$, for $0<i < L$. 
By Corollary \ref{lcusps}, this is is asymptotically
$$e^{-\sum_{i=0}^{L-1} {{1}\over{i}}} \sim {{e^{-\gamma}}(L-1)^{-1}},$$
 where $\gamma$ is 
Euler's constant.

Hence the lemma is proved if we replace the conclusion ``probability $\to
1$'' with the weaker conclusion ``probability $\ge (\const) > 0\
{\hbox{as}}\ n \to \infty$.'' This is a sufficiently strong version of the
lemma to obtain the results announced in \cite{BM}.

We now show how to obtain the stronger results of Theorem
\ref{cusps}A. To that end, suppose that the cusp $C_0$ of the 
surface $S^O(\Gamma,
\calO)$ has a canonical horocycle of length $< L$. We would like to
choose a larger horocycle about this cusp. 

There are two obstructions to choosing such a larger horocycle. The
first obstruction is that as we increase the length of the horocycle
about $C_0$, it may cease to be injective. This will happen if there
is a short closed geodesic in the neighborhood of $C_0$. 

The second obstruction is that, as we increase the length of the
canonical horocycle about $C_0$, we must decrease the lengths of the
horocycles of nearby cusps in order to keep the interiors of the
horocycles disjoint. When we decrease the length of a nearby
horocycle, it may then cease to have length $ \ge L$.

Both of these considerations are handled by Lemma \ref{dist}. Indeed,
both of these obstructions arise from the possibility that in the
graph $\Gamma$, there may be a short closed path close to the left-hand-turn
path corresponding to $C_0$. According to Lemma \ref{dist}, the
probability of this occurring is asymptotically $0$.

This argument is illustrated in Figure \ref{horo} below.
The cusp in question lies between the vertical lines $x=0$
and $x=2$, and the canonical horocycle, the line $y=1$, has
length $2$. We increase its length to $8$ by lowering the
horocycle to the line $y=1/4$. This will be possible if none
of the points $x+iy$ with $0 < x<2, y> 1/4$ are identified
in the surface (the first obstruction), and if all the
horocycles which meet the line $y= 1/4$ have images in the
surface which are sufficiently long (the second obstruction).

This concludes the proof of the lemma.

\begin{figure}[!h]
\begin{center}

\includegraphics[angle=0,scale=0.6]{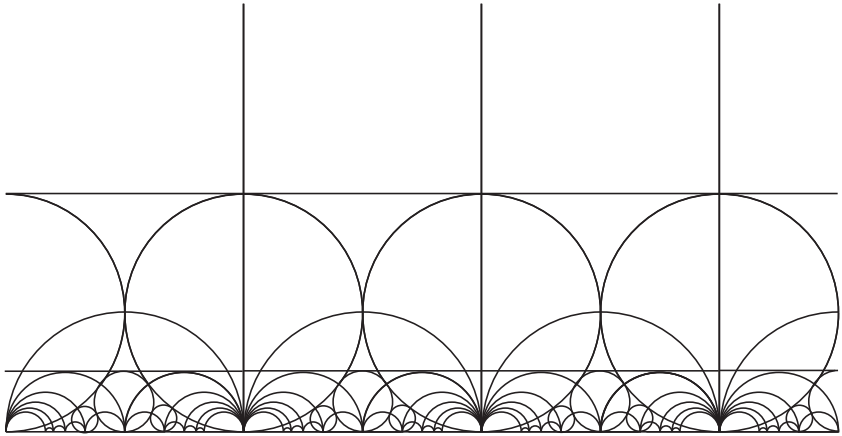}

\caption{Increasing the size of horocycle}\label{horo}
\end{center}
\end{figure}

\ifx\undefined\bysame
\newcommand{\bysame}{\leavevmode\hbox to3em{\hrulefill}\,}
\fi

\end{document}